\documentclass{amsart}
\usepackage{amsmath}
\usepackage{amsthm}
\usepackage{amssymb}
\usepackage{latexsym}               % LaTeX Symbol Font,amssymb
\usepackage{amsgen}
\usepackage[mathscr]{eucal}
\numberwithin{equation}{section}
\usepackage[dvipdfm, usenames]{color}

\newcommand{\al}{\alpha}

\newcommand{\e}{\epsilon}

\newcommand{\ti}{\tilde}

\makeatother

\newtheorem{lem}{Lemma}[section]
\newtheorem{thm}[lem]{Theorem}

\begin{document}

\title
{Inviscid Limit Problem of radially symmetric stationary solutions for compressible Navier-Stokes equation}
\author{Itsuko Hashimoto$\dagger$}
\email{itsuko@se.kanazawa-u.ac.jp}
\address{Kanazawa university, Osaka city university, Japan}

\author{Akitaka Matsumura}
\email{akitaka@math.sci.osaka-u.ac.jp}
\address{Osaka university, Japan}

\thanks{
$\dagger$ \ Supported by JSPS Grant No. 17K14227
}

\date{}
\maketitle

\def\cal#1{{\fam2#1}}

%%%%%%    TEXT START    %%%%%%
%\begin{center}
%\small{\it Dedicated to Professor Ling Hsiao on the occasion of her 80th birthday}
%\end{center}

%%%%%%%%%%%%%%%%%%%%%%%%  abstract  %%%%%%%%%%%%%%%%%%%%%%%%%%%%%

\begin{abstract} 
The present paper is concerned with an inviscid limit problem of 
radially symmetric stationary solutions for an exterior problem in  
$\mathbb{R}^n (n\ge 2)$ to  compressible Navier-Stokes equation, 
describing the motion of viscous barotropic gas without external forces, 
where boundary and far field data are prescribed. For both inflow and outflow problems, 
the inviscid limit is considered in a suitably small neighborhood of the far field state.
For the outflow problem, we prove the uniform convergence 
of the Navier-Stokes flow toward the corresponding Euler flow in the inviscid limit. 
On the other hand, for the inflow problem, we show that the Navier-Stokes flow 
uniformly converges toward a linear superposition of 
the corresponding boundary layer profile and the
Euler flow in the inviscid limit. 
The estimates of algebraic rate toward the inviscid limit are also obtained.
\end{abstract}

\medskip

%%%%%%%%%%%%%%%%%%%%%%%  Section 1  %%%%%%%%%%%%%%%%%%%%%%%%%%%%%%
%%%%%%%%%%%%%%%%%%%%%%%  introduction  %%%%%%%%%%%%%%%%%%%%%%%%%%

\section{Introduction and Main theorem}
\noindent
In the present paper, we consider an inviscid limit problem of radially symmetric stationary solutions of the compressible Navier-Stokes equation which describes a barotropic motion of viscous gas in the exterior domain 
$\Omega$ to a ball in $\mathbb{R}^{n}\,(n\ge 2)$:
\begin{eqnarray}
 \label{cns} 
 \left\{
\begin{array}{l}
\rho_t+{\rm div} (\rho U)=0, \\[2mm]
(\rho U)_t+{\rm div} (\rho U\otimes U)+\nabla p =
\nu \bigtriangleup U +(\nu+\lambda)\nabla({\rm div} U),\quad
t>0,\ x \in \Omega,
\end{array}
 \right.\,
\end{eqnarray}
where $\Omega=\{x \in\mathbb{R}^{n}\,(n\ge 2); |x| > r_0\}$\,($r_0$ is a positive constant),
$\rho=\rho(t,x)>0$ is the mass density, $U=(u_1(t,x),\cdots ,u_n(t,x))$ is the fluid velocity, 
and $p=p(\rho)$ is the pressure given by a smooth function of 
$\rho$ satisfying $p'(\rho)>0\ (\rho>0)$.
Furthermore, $\nu$ and $\lambda$ are 
the shear and second viscosity coefficients respectively, which are assumed to be 
constants satisfying $\nu>0, 2\nu+ n\lambda \ge 0$. 
Our concern is the problem of which solution is given by the radially symmetric solution: 
\begin{align}
\label{t}
\rho(t,x)=\rho(t,r), \quad U(t,x)=\frac{x}{r}\,u(t, r),\quad r= |x|,
\end{align}
where $\rho(t, r)$ and  $u(t, r)$ are scalar functions. 
Under the assumption (\ref{t}), the equation 
for the radially symmetric solution $(\rho, u)$ of (1.1) becomes 
\begin{eqnarray}
 \label{nr} 
 \left\{\begin{array}{l}
%   \rho_{t}+(\rho u)_r+(n-1)\dfrac{\rho u}{r}=0, \\
(r^{n-1}\rho)_{t} + {(r^{n-1}{\rho}{u})_r}= 0, \\
  (\rho u)_t+(\rho u^2+p(\rho))_r+(n-1)
\dfrac{\rho u^2}{r}=\mu \big(\dfrac{(r^{n-1}{u})_r}{r^{n-1}}\big)_r,\ \ \ t>0,\ r>r_0,
 \end{array}
 \right.
\end{eqnarray}
where $\mu = 2\nu + \lambda>0$. Now, we consider the initial boundary value problems
to (\ref{nr}) under the initial condition
% is prescribed to be spatial asymptotically constant
\begin{equation}
(\rho, u)(0,r)=(\rho_0, u_0)(r),\quad r >r_0, 
\end{equation}
the far field condition
\begin{equation}
\displaystyle{\lim_{r \to \infty}}(\rho, u)(t,r)=(\rho_+,u_+),\quad t>0,
\end{equation}
% \begin{equation}
% (\rho, u)(0,r)=(\rho_0, u_0)(r),\ r >r_0, \quad
% \end{equation}
and also the following two types of boundary conditions depending on the sign of the velocity
on the boundary
\begin{equation}
\label{bc}
 \left\{
\begin{array}{ll}
(\rho, u)(t,r_0) =(\rho_-, u_-),\quad t>0,& (u_->0),\\[10pt]
u(t,r_0) = u_-,\quad t>0,& (u_-\le 0),
\end{array}
 \right.
\end{equation}
where $\rho_{\pm}>0, u_{\pm}$ are given constants. 
The case $u_->0$ is known as ``inflow problem'', the case $u_-=0$ as 
``impermeable wall problem'', and the case $u_-<0$ as ``outflow problem''. 
In the previous paper \cite{H-M}, for both inflow and outflow problems, 
we showed the existence of a unique radially symmetric stationary solution 
for (\ref{nr})-(\ref{bc}) in a suitably small neighborhood of the far field state. 
For the outflow problem, we showed in \cite{NSH} the asymptotic stability 
of stationary solutions for (\ref{nr})-(\ref{bc}) 
provided  the initial disturbance is  small. 
The stationary solution of the equation (\ref{nr}) is independent of time $t$ 
and satisfies the same boundary condition (1.5) and (1.6). 
Therefore, the stationary solution verifies 
\begin{align}
 \label{nrs0} 
 \left\{
\begin{aligned}
  &\ (r^{n-1}\rho u)_r = 0,   \\
  &\ \rho uu_r+p(\rho)_r=\mu(\frac{(r^{n-1}u)_r}{r^{n-1}})_r, \qquad r>r_0,\\[2mm]
  &\ \displaystyle{\lim_{r\to \infty}} (\rho, u)(r) =(\rho_+, u_+), \\[2mm]
  &\ (\rho, u)(r_0)=(\rho_-, u_-)\ \ (u_->0), \quad u(r_0)=u_- \ \ (u_- \le 0).
\end{aligned}
 \right.\,
\end{align}
From the first equation 
in (\ref{nrs0}), we easily see it holds
\begin{align}
\label{n1}
r^{n-1}\rho(r) u(r)=\epsilon, \qquad r\ge r_0, 
%\qquad and \qquad  r_0^{n-1}\rho_- u_-=\epsilon,
\end{align}
for some constant $\epsilon$, and it also holds from the boundary conditions that
\begin{equation}
\label{Nep}
\epsilon = r^{n-1}_0\rho_- u_-\quad(u_->0),\qquad 
\epsilon = r^{n-1}_0\rho (r_0) u_-\quad(u_-\le0).
\end{equation}
The formula (\ref{n1}) implies that if $n \ge 2$, 
\begin{align*}
u_+=\lim_{r\to \infty}u(r)= \lim_{r\to \infty}\frac{\epsilon}{r^{n-1}\rho_+}=0.
\end{align*}
Hence, we need to assume $u_+=0$ in the present paper. 
Here we recall the results in the previous papers 
\cite{H-M} and \cite{NSH} which gives the existence of the solution of (\ref{nrs0}).
%%%%%%%%%%%%%%%%%%%%%%%%%%%%%%%%%%%%%%%%%%%%%%%%%%%%%%%%%%%%%%%%%5
\bigskip
\noindent
\begin{thm}\rm{(\cite{H-M,NSH})}
\label{mt}
\quad
Let $n \ge 2$ and $u_+=0$. 
Then, for any $\rho_+>0$,
there exist positive constants $\epsilon_0$ and $C$ satisfying the following:

\medskip

\noindent
{\rm (I)} \ Let $u_->0$. If $|u_-|+|\rho_--\rho_+| \le \epsilon_0$, 
there exists a unique smooth solution $(\rho,u)$ of 
the problem {\rm (\ref{nrs0})} satisfying 
\begin{align*}
\begin{aligned}
  & |\rho(r)-\rho_+|\le C r^{-2(n-1)}(|u_-|^2+|\rho_--\rho_+|),\\[5pt]
  & C^{-1}r^{-(n-1)}|u_-| \le |u(r)| \le C r^{-(n-1)}|u_-|, \qquad r\ge r_0.
\end{aligned}
\end{align*}\\
{\rm (II)} \ Let $u_-\le 0$. If $|u_-|\le \epsilon_0$, 
there exists a unique smooth solution $(\rho,u)$ 
of the problem {\rm (\ref{nrs0})} satisfying
\begin{align*}
\begin{aligned}
  & |\rho(r)-\rho_+|\le C r^{-2(n-1)}|u_-|^2, \\[5pt]
  & C^{-1}r^{-(n-1)}|u_-| \le |u(r)|\le C r^{-(n-1)}|u_-|, \qquad r\ge r_0.
\end{aligned}
\end{align*}
%where $C$ is a constant independent of $\mu$.
\end{thm}
\noindent
Here it is noted that the constants $\epsilon_0$ and $C$ in the above statements 
may depend on $\mu$. However, in order to consider the inviscid limit problem,
we can show the above statements uniformly hold for $\mu \in (0,1]$ (see Theorem 2.1).

 In the present paper, we shall consider the inviscid limit problem  
 of radially symmetric stationary solutions 
 for the compressible Navier-Stokes equation, that is, 
 clarify the relation of the solution 
 of (1.7) to that of the Euler equation:
% 
% 
%To be more specific, for outflow problem, we show that the solution to the problem (\ref{nrs0}) converges to the corresponding stationary solution for the Euler equation:
\begin{align}
 \label{Es} 
 \left\{
\begin{aligned}
  &\ (r^{n-1}\rho^{E} u^{E})_r = 0,   \\[2mm]
  &\ \rho^{E} u^{E}u^{E}_r+p(\rho^{E})_r=0, \qquad r>r_0,\\[2mm]
  &\ \displaystyle{\lim_{r\to \infty}} (\rho^{E}, u^{E})(r) =(\rho_+, u_+), 
\end{aligned}
 \right.\,
\end{align}
with the corresponding suitable boundary condition for
$(\rho^{E},u^{E})(r_0)$.
%in a suitably small neighborhood of the far field state. 
%On the other hand, for inflow problem, we show 
%that the solution to the problem (\ref{nrs0}) 
%converges to the superposition of corresponding 
%stationary solution for the Euler equation and boundary layer solution. 
Here it is noted that when we consider the initial boundary value 
problem for the Euler equation, as long as the velocity $u$ is small,
that is, the system is subsonic, only single boundary condition 
can be allowed for the problem to be well-posed.
In this sense, for the outflow problem, we naturally 
can expect the corresponding
boundary condition $(\rho^{E},u^{E})(r_0)$ should be $u^{E}(r_0)=u_-$, and,
in fact, we shall show that  
the solution for the problem (\ref{nrs0}) uniformly converges 
toward that for the Euler equation (\ref{Es}) with the boundary condition
$u^{E}(r_0)=u_-$ as $\mu \to 0$.
On the other hand, for the inflow problem, no mater how we choose
the boundary condition for the Euler equation, because of  the inconsistency
of the number of the boundary conditions, we expect  a boundary layer 
does appear. First, we naturally choose the single boundary condition 
for the Euler equation as $(\rho^{E}u^{E})(r_0)= \rho_-u_-$, and next, 
introduce a boundary layer solution $(\rho^B, u^B)$
 for  (\ref{nrs0}) to smoothly connect the 
boundary data $(\rho_-,u_-)$ and the boundary value 
 $(\rho^{E},u^{E})(r_0)$ of the
solution of the
Euler equation.
To do that, taking scale transformation of space as
 $y:=\frac{r-r_0}{\mu}$, rewriting (\ref{nrs0}) in terms of 
 $y$, and formally taking $\mu=0$, we can deduce the
 following boundary layer equation (see Section 4.2 for details):
 %of space as $y:=\frac{r-r_0}{\mu}$ and taking inviscid limit. 
%The boundary layer equation to the equation (\ref{nrs0}) is 
%given in the case $u_->0$ and derived by scale transformation 
%of space as $y:=\frac{r-r_0}{\mu}$ and taking inviscid limit. 
%Therefore the boundary layer solution $(\rho^B, u^B)$ verifies 
\begin{align}
 \label{Ebl} 
 \left\{
\begin{aligned}
  &\ (\rho^Bu^B)_y = 0,   \\
  &\ \rho^Bu^Bu^B_y+p(\rho^B)_y=u^B_{yy}, \qquad y>0,\\[2mm]
  &\ \displaystyle{\lim_{y\to \infty}} (\rho^B,u^B)(y) =(\rho^{E}(r_0), u^{E}(r_0)), \\[2mm]
  & (\rho^B,u^B)(0)=(\rho_-, u_-).
\end{aligned}
 \right.\,
\end{align}
%where $r=1+\mu y$.
Then we can show 
that the solution to the problem (\ref{nrs0}) 
uniformly converges toward a linear superposition of the corresponding 
solution of the Euler equation and the boundary layer solution.

There have been many works on the inviscid limit problems for 
the Navier-Stokes equations, for example, Hoff-Liu \cite{HL}, Williams \cite{WW}, Sueur \cite{S} in the compressible cases, and Constantin-Chemin \cite{C,J-PC}, 
Asano-Sammartino \cite{A,S-I,S-II}, 
Maekawa-Mazzucato \cite{MM} in 
the incompressible cases. 
However, it seems there have been no results on the radially 
symmetric solutions of the inflow/outflow 
problems to the compressible Navier-Stokes equation. 
In particular, for the inflow problem, we show 
a boundary layer for both density and velocity 
does appear in the inviscid limit, which phenomenon 
is not observed in the previous papers. 

Now we are ready to state the main results in the present paper. The first theorem is concerned  with the inviscid limit  
problem  for the outflow problem. 
%of the existence of the stationary solution for the boundary value problem of E%uler equation (\ref{Es}).
%%%%%%%%%%%%%%%%%%%%%%%%%%%%%%%%%%%%%%%%%%%%%%%%%%%%%%%%%%%%%%%%%5
%\bigskip
%\noindent
%\begin{thm}\cite{HMn1}
%\label{mt}
%Let $n \ge 2$ and $u_+=0$. 
%Then, for any $\rho_+>0$,
%there exist positive constants $\epsilon_0$ and $C$ satisfying the following:
%
%\medskip
%
%\noindent
%{\rm (I)} \ Let $u_->0$. If $|u_-|+|\rho_--\rho_+| \le \epsilon_0$, 
%there exists a unique smooth solution $(\rho,u)$ of 
%the problem {\rm (\ref{nrs0})} satisfying 
%\begin{align}
%\label{1}
%\begin{aligned}
%  & |\rho(r)-\rho_+|\le Cr^{-(n-1)}(|u_-|^2+|\rho_--\rho_+|),\\[5pt]
%  & C^{-1}r^{-(n-1)}|u_-| \le |u(r)| \le Cr^{-(n-1)}|u_-|, \qquad r\ge r_0.
%\end{aligned}
%\end{align}
%Furthermore, for any positive constant $h$, there exists a positive constant 
%$C_h$ such that it holds
%\begin{align}
%\label{bl}
%\sup_{r\ge r_0+h}|\rho(r)-\rho_+| \le C_h|u_-|^2.
%\end{align} \\
%{\rm (II)} \ Let $u_-\le 0$. If $|u_-|\le \epsilon_0$, 
%there exists a unique smooth solution $(\rho,u)$ 
%of the problem {\rm (\ref{nrs0})} satisfying
%\begin{align}
%\begin{aligned}
%  & |\rho(r)-\rho_+|\le Cr^{-2(n-1)}|u_-|^2, \\[5pt]
%  & C^{-1}r^{-(n-1)}|u_-| \le |u(r)|\le Cr^{-(n-1)}|u_-|, \qquad r\ge r_0.
%\end{aligned}
%\end{align}
%\end{thm}
\begin{thm}
\label{etE}
Let $n \ge 2$, $u_- \le 0$ and $u_+=0$. Then, for any $\rho_+>0$, 
there exist positive constants $\epsilon_0$ and $C$ 
which are independent of $\mu \in (0,1]$ satisfying the following:

\medskip

\noindent
{\rm (I)} \ If  $|u_-|\le \epsilon_0$, there exists a unique smooth solution 
$(\rho^{E},u^{E})$ of the problem {\rm (\ref{Es})} with the boundary condition $u^{E}(r_0)=u_-$ satisfying
\begin{align*}
\begin{aligned}
  & |\rho^{E}(r)-\rho_+|\le Cr^{-2(n-1)}|u_-|^2, \\[5pt]
  & C^{-1}r^{-(n-1)}|u_-| \le |u^{E}(r)|\le Cr^{-(n-1)}|u_-|, \qquad r\ge r_0.
\end{aligned}
\end{align*}\\
{\rm (II)} \ If $|u_-|\le \epsilon_0$, the solution $(\rho,u)$ of 
the problem {\rm (\ref{nrs0})} uniformly converges 
toward the solution $(\rho^{E},u^{E})$ of the problem {\rm (\ref{Es})} 
as $\mu \to 0$, that is,
it holds the estimate
\begin{align*}
\begin{aligned}
%\sup_{1\le r}|\rho(r)-\rho^{E}(r)| \to 0, \quad   \sup_{1\le r}|u(r)-u^{E}(r)| \to 0, \quad (\mu \to 0).
\sup_{r\ge r_0}|(\rho-\rho^{E}, u-u^{E})(r)| \le C\mu.
\end{aligned}
\end{align*}

\end{thm}
The second theorem is with the inviscid limit  problem  for the inflow problem. 
\begin{thm}
\label{bl+E}
Let $n \ge 2$, $u_->0$ and $u_+=0$. Then, for any $\rho_+>0$,
there exist positive constants $\epsilon_0$, $\nu_0$, and $C$ which are independent of $\mu \in (0,1]$ satisfying the following:

\medskip

\noindent
{\rm (I)} \ If  $|u_-|+|\rho_--\rho_+| \le \epsilon_0$, there exists a unique smooth solution $(\rho^{E},u^{E})$ of the problem {\rm (\ref{Es})} with the boundary condition $u(r_0)\rho(r_0)={u_-\rho_-}$ satisfying
\begin{align*}
\begin{aligned}
  & |\rho^{E}(r)-\rho_+|\le Cr^{-2(n-1)}|u_-|^2, \\[5pt]
  & C^{-1}r^{-(n-1)}|u_-| \le |u^{E}(r)|\le Cr^{-(n-1)}|u_-|, \qquad r\ge r_0.
\end{aligned}
\end{align*}

\noindent
{\rm (II)} \ If  $|\rho_--\rho_+|+|u_-|\le \epsilon_0$, 
there exists a unique smooth solution $(\rho^B,u^B)$ 
of the problem {\rm (\ref{Ebl})} satisfying
\begin{align*}
\begin{aligned}
  & |{\rho}^B(y)-\rho^{E}(r_0)|
  \le Ce^{-\frac{\nu_0}{\epsilon} y}(|u_-|^2+|\rho_--\rho_+|), \\[5pt]
  & |{u}^B(y)-u^{E}(r_0)|\le Ce^{-\frac{\nu_0}{\epsilon} y}(|u_-|^2+|\rho_--\rho_+|),
  \qquad y\ge 0,
\end{aligned}
\end{align*}
where $\epsilon=\rho_-u_-$.

\noindent
{\rm (III)} \ 
If $|\rho_--\rho_+|+|u_-|\le \epsilon_0$, the solution $(\rho,u)$ of 
the problem {\rm (\ref{nrs0})} uniformly converges 
toward a linear superposition of the boundary layer solution $(\rho^B,u^B)$ 
of the problem {\rm (\ref{Ebl})} and the solution $(\rho^{E},u^{E})$ 
of the problem {\rm (\ref{Es})} as $\mu \to 0$, that is,
it holds the estimates for $\al \in (0,1)$
\begin{align*}
\begin{aligned}
&\sup_{r \ge r_0}|\rho(r)-\rho^{E}(r)-\rho^B(\dfrac{r-r_0}{\mu})+\rho^{E}(r_0)|
\le C\mu^\al,\\
&\sup_{r \ge r_0}|u(r)-u^{E}(r)-u^B(\dfrac{r-r_0}{\mu})+u^{E}(r_0)|
\le C\mu^\al,
\end{aligned}
\end{align*}
where $C$ also depends on $\alpha\in (0,1)$.
\end{thm}
This paper is organized as follows. 
In Section 2, we reformulate the Navier-Stokes equation (\ref{nrs0}) 
and state the existence of the solution and its decay rate 
estimates which are uniform with respect to  $\mu$. In Section 3,
we consider the outflow problem and show that 
the inviscid limit of the solution of (\ref{nrs0}) is 
given by the solution of the Euler equation  (\ref{Es}).  
In Section 4, we finally prove that 
for the inflow problem,
the inviscid limit of the solution of (\ref{nrs0}) is given by
a linear superposition of the solution of Euler equation (\ref{Es}) 
and the boundary layer solution of (\ref{Ebl}).

\section{Reformulation of Navier-Stokes equation}
In this section, we reformulate the problem (\ref{nrs0}) as in the previous paper \cite{H-M}.  
First, we assume $u_- \ne 0$ in what follows, because 
if $u_-= 0$, that is, $\epsilon =0$ in (\ref{Nep}), 
the unique solution of (\ref{nrs0}) easily turns out to be the trivial
constant state $(\rho,u) \equiv (\rho_+,0)$. 
We also assume $r_0=1$ without loss of generality. 
Next, introduce the specific volume $v$ by $v=1/\rho$ 
(accordingly, denote $v_\pm$ by $1/\rho_\pm$). Then, by (\ref{n1}),
the velocity $u$ is given in terms of $v$ as
\begin{align}
\label{urv}
u(r)=\frac{\epsilon}{r^{n-1}}v(r),\quad r\ge 1,
\end{align}
where $\epsilon = u_-/v_-\ (u_->0)$, and
$\epsilon = u_-/v(1)\ (u_- <0)$.
Substituting (2.1) into the second equation of (\ref{nrs0}), we have
\begin{align}
\label{base}
\frac{\epsilon^2}{r^{n-1}}\left(\frac{v}{r^{n-1}}\right)_r+\tilde{p}(v)_r
= \epsilon\mu\left(\frac{v_r}{r^{n-1}}\right)_r,
\end{align}
where $\tilde{p}(v) := p(1/v)$, and it holds $\tilde{p}'(v)<0\,(v>0)$ by the assumption on $p(\rho)$.
Now, we further introduce a new unknown function $\eta$, as the deviation of $v$ from the
far field state $v_+$, by
\begin{equation}
\label{eta}
\eta(r) = v(r) - v_+,\qquad r \ge 1.
\end{equation}
Plugging (\ref{eta}) into (\ref{base}), we obtain 
\begin{equation}
\label{Narrange}
\epsilon \mu \left(\frac{ \eta_r}{r^{n-1}}\right)_r=
\tilde p(v_++\eta)_r
+\frac{\epsilon^2 v_+}{2}\Big(\frac{1}{r^{2(n-1)}}\Big)_r
+\frac{\epsilon^2}{r^{n-1}}\left(\frac{\eta}{r^{n-1}}\right)_r,
\end{equation}
where $\epsilon = u_-/v_-\ (u_->0)$, and $\epsilon = u_-/(v_++\eta(1))\ (u_- < 0)$.
Under the far field condition $\eta(\infty)=0$, 
the equation (\ref{Narrange}) is also equivalent to
\begin{equation}
\label{Narranger}
\left(\epsilon \mu \frac{\eta_r}{r^{n-1}}-\tilde p(v_++\eta) 
-\frac{\epsilon^2 v_+}{2r^{2(n-1)}}
-\frac{\epsilon^2\eta}{r^{2(n-1)}}
+\epsilon^2(n-1)\int_r^{\infty}\frac{\eta(s)}{s^{2n-1}}ds
\right)_r =0,
\end{equation}
which implies that the function in the parenthesis in the
left hand side of (\ref{Narranger}) is identically equals to a constant
$c_0$ for $r >1$.
Then, it follows from the far field condition that $c_0=-\tilde p(v_+)$. Thus, we finally 
have the following reformulated problem in terms of
$\eta$: 
\begin{align}
\label{reform}
 \left\{
\begin{aligned}
  &\eta_r= 
\frac{r^{n-1}}{\epsilon\mu}\big(\tilde p(v_++\eta)- \tilde p(v_+)\big)\\[2mm]
&\qquad\ +\frac{\epsilon v_+}{2\mu}\frac{1}{r^{n-1}}
+\frac{\epsilon\eta}{\mu r^{n-1}}
-\frac{\epsilon(n-1)r^{n-1}}{\mu}\int_r^{\infty}\frac{\eta(s)}{s^{2n-1}}ds,
\quad r>1,\\[2mm]
 &\displaystyle{\lim_{r\to \infty}}\eta(r) =0,  \\[2mm]
  &\eta(1)= \eta_-:=v_--v_+\ (u_->0),\quad no\ boundary\ condition\ (u_- < 0),
 \end{aligned}
 \right.
\end{align}
where $\epsilon = u_-/v_-\ (u_->0)$, and $\epsilon = u_-/(v_++\eta(1))\ (u_- < 0)$.
Once the desired solution $\eta$ of (\ref{reform}) is obtained, 
the velocity $u$ is immediately obtained by (\ref{urv}) as
\begin{equation*}
u(r)=\frac{u_-(v_++\eta(r))}{v_-r^{n-1}}\quad (u_->0),
\quad
u(r)=\frac{u_-(v_++\eta(r))}{(v_++\eta(1))r^{n-1}}\quad (u_- < 0).
\end{equation*} 
The theorem for the reformulated problem (\ref{reform}) given in the previous papers is as follows,
where we newly added the statement on the dependency of constants on $\mu$.
\noindent
\begin{thm}{\rm(\cite{H-M,NSH})}
\label{rmt}
\quad Let $n \ge 2$. Then, for any $v_+>0$, there exist positive constants $\epsilon_0$ and $C$ which are
 independent of $\mu \in (0,1]$ satisfying the following:

\medskip

\noindent
{\rm (I)} \ Let $u_->0$. If $|u_-|+|\eta_-|\le \epsilon_0$, 
there exists a unique smooth solution $\eta$ of the problem {\rm (\ref{reform})} satisfying 
\begin{align*}
|\eta(r)|\le C r^{-(n-1)}(|u_-|^2+|\eta_-|),\qquad  r\ge 1. 
\end{align*}
%Furthermore, for any positive constant $h$, 
%there exists a positive constant $C_h$ satisfying
%\begin{align*}
%\sup_{r\ge r_0+h}|\eta(r)| \le C_h|u_-|^2.
%\end{align*} \\
{\rm (II)} \ Let $u_- < 0$. If $|u_-| \le \epsilon_0$, 
there exists a unique smooth solution $\eta$ of the problem {\rm (\ref{reform})} satisfying
\begin{align}
\label{estouteta}
|\eta(r)|\le C r^{-2(n-1)}|u_-|^2, \qquad r\ge 1. 
\end{align}
\end{thm}
\noindent
%We can easily see that the Theorem \ref{mt} is a direct  consequence 
%by Theorem \ref{rmt}.
As for the proof of Theorem 2.1, since the argument is entirely same as in
\cite{H-M,NSH} except keeping in mind  the dependency on $\mu$ 
a bit carefully,
so we omit it.
%%%%%%%%%%%%%%%%%%%%%%%%%%%%%%%%%%%%%%%%%%%%%%%%%%%%%%%%

\section{Outflow problem}
In this section, we consider the case $u_-<0$ and 
show the statement (II) in Theorem \ref{etE}. That is, we show that 
the inviscid limit of the solution of Navier-Stokes equation (\ref{nrs0}) is 
given by
that of Euler equation (\ref{Es}) with boundary condition $u^{E}(r_0)=u_-$. 
%%%%%%%%%%%%%%%%%%%%%%%%%%%%%%%%%%%%%%%%%%%%%%%%%%%%%%%%
\subsection{Euler equation} 
In this subsection, we reformulate the problem (\ref{Es}) 
with the  boundary condition $u^{E}(r_0)=u_-$,
and show the existence of the solution. From the first equation 
of (\ref{Es}), we easily see that it holds
\begin{align*}
%\label{n1e}
r^{n-1}\rho^{E}(r) u^{E}(r)=\epsilon^{E}, \qquad r\ge 1, 
%\qquad and \qquad  r_0^{n-1}\rho_- u_-=\epsilon,
\end{align*}
for some constant $\epsilon^{E}$, and it also holds from the boundary condition that
\begin{equation}
\label{Eep}
\epsilon^{E} =  u_- \rho^{E}(1)= u_-/v^{E}(1),
\end{equation}
where $v^{E}(r):=1/\rho^{E}(r)$.
On the other hand, by the first equation in (\ref{Es}), 
the velocity $u^{E}$ is given by
\begin{align}
\label{uE}
u^{E}(r)=\frac{\epsilon^{E}}{r^{n-1}}v^{E}(r),\quad r\ge 1.
\end{align}
Substituting (\ref{uE}) into the second equation of (\ref{Es}) and introducing a new variable $\eta^{E}(r):=v^{E}(r)-v_+$, we have
\begin{align}
\label{baseE}
\frac{|\epsilon^{E}|^2}{r^{n-1}}\left(\frac{v_++\eta^{E}}{r^{n-1}}\right)_r+\tilde{p}(v_++\eta^{E})_r=0, 
\end{align}
which also can be rewritten as
\begin{align}
\label{E1}
\eta^{E}_r(r)=\frac{|\epsilon^{E}|^2(n-1)(v_++\eta^{E})}{(|\epsilon^{E}|^2r^{-2(n-1)}+\tilde p'(v_++\eta^{E}))r^{2n-1}}.
\end{align}
Integrating (\ref{E1}) in terms of space variable $r$ on $[r,\infty)$, we obtain 
\begin{align}
\label{E2}
\begin{aligned}
\eta^{E}(r)&=\int_r^{\infty}\frac{|\epsilon^{E}|^2(n-1)(v_++\eta^{E})}{(|\tilde p'(v_++\eta^{E})|-|\epsilon^{E}|^2s^{-2(n-1)})s^{2n-1}}ds\\
&:=\int_r^{\infty}G[\eta^{E}](s)ds,
\end{aligned}
\end{align}
%\begin{align}
%\label{E2}
%\eta^{E}(r)=-\int_r^{\infty}\frac{|\epsilon^{E}|^2(n-1)(v_++\eta^{E})}{(|\epsilon^{E}|^2s^{-2(n-1)}+\tilde p'(v_++\eta^{E}))s^{2n-1}}ds,
%\end{align}
where we used the far field condition of $\displaystyle{\lim_{r\to \infty}}\eta^{E}(r)=0$. On the other hand, applying the same strategy as (\ref{base}) to (\ref{reform}), we also can reformulate (\ref{Es}) to the following problem:
\begin{align}
\label{reformE}
 \left\{
\begin{aligned}
&\big(\tilde p(v_++\eta^{E})- \tilde p(v_+)\big)+\frac{|\epsilon^{E}|^2 v_+}{2}\frac{1}{r^{2(n-1)}}\\[2mm]
&\qquad\qquad\qquad\qquad\quad
+\frac{|\epsilon^{E}|^2\eta^{E}}{r^{2(n-1)}}-|\epsilon^{E}|^2(n-1)\int_r^{\infty}\frac{\eta^{E}(s)}{s^{2n-1}}ds=0,
\quad r \ge 1,
\\[2mm]
&\displaystyle{\lim_{r\to \infty}}\eta^{E}(r) =0,& 
\end{aligned}
 \right.
\end{align}
where $\epsilon^{E} = u_-/(v_++\eta^{E}(1))$. 
The existence theorem for the reformulated problem (\ref{E2})/(\ref{reformE}) which we need to prove is
\begin{thm} 
\label{rmtE}
Let $n \ge 2$. Then, for any $v_+>0$, there exist 
positive constants $\epsilon_0$ and $C$ such that if $|u_-| \le \epsilon_0$, 
there exists a unique smooth solution $\eta^{E}$ of the problem {\rm (\ref{E2})/(\ref{reformE})} satisfying
\begin{subequations}
\begin{align}
\label{decayE1}
&|\eta^{E}(r)|\le Cr^{-2(n-1)}|u_-|^2, \qquad r\ge 1, \\[2mm]
\label{decayE2}
&|\eta^{E}_r(r)|\le Cr^{-2n+1}|u_-|^2, \ \ \qquad r\ge 1.
\end{align}
\end{subequations}
\end{thm}
Once Theorem 3.1 is proved, since the velocity $u^{E}$ is given by (3.2),
thus the statement (I) in Theorem 1.2 is also proved.

To prove the existence of the solution with the decay rate estimate (\ref{decayE1}), 
we look for a solution of (\ref{E2}) in the Banach space $X$, 
with its norm $\|\cdot\|_X$, defined by
\begin{align*}
%\label{eX}
X=\{\eta^{E} \in C([1,\infty));\ \sup_{r\ge 1}|r^{2(n-1)}\eta^{E}(r)|< \infty \},
\quad \|\eta^{E}\|_X = \sup_{r\ge 1}|r^{2(n-1)}\eta^{E}(r)|.
\end{align*}
\noindent
% By using (\ref{eq4}), we make the sequence of iteration as
To do that, we construct the approximate sequence $\{{\eta^{E}}^{(m)}\}_{m\ge 0}$ by
\begin{align}
 \label{eas} 
 \left\{\begin{aligned}
&{\eta^{E}}^{(0)}(r)=\int_r^{\infty}G[0](s) ds,\quad  r\ge 1, \\[2mm]
&{\eta^{E}}^{(m+1)}(r)=\int_r^{\infty} G[{\eta^{E}}^{(m)}](s) ds,
\ \ (m \ge 0),  \quad r\ge 1.\\[2mm]
 \end{aligned}
 \right.\,
\end{align}

{\it Proof of Theorem \ref{rmtE}}
\quad We show that $\{{\eta^{E}}^{(m)}\}_{m\ge 0}$ is a Cauchy sequence in $X$
 for suitably small $|u_-|$.
To do that,  we first show the uniform boundedness of 
${\eta^{E}}^{(m)}\ (m\ge 0)$ in $X$ for suitably small 
$|u_-|$. More precisely, we show that for any fixed $v_+$,
there exist positive constants $\epsilon_0$ and $C$  
such that if $|u_-|\le \epsilon_0$, there exists a positive constant $M$ satisfying
\begin{equation}
\label{eu}
\|{\eta^{E}}^{(m)}\|_X \le M \le C|u_-|^2,\quad m\ge 0.
\end{equation}
Here and in what follows, we use the letter $C$ and $\epsilon_0$ to denote generic
positive constants which are independent of $u_-$, 
but may depend on $v_+$ and $n$. For the proof, we first assume
\begin{equation}
\label{econdi0}
\frac{|\tilde p'(\frac{3}{2}v_+)|}{2} \le |\tilde p'(\frac{3}{2}v_+)|-|\frac{2u_-}{v_+}|^2,
\end{equation}
which is equivalent to
\begin{equation}
\label{econdi1}
|u_-| \le \frac{1}{2}\sqrt{\frac{|\tilde p'(\frac{3}{2}v_+)|}{2}}v_+.
\end{equation}
Let us show (\ref{eu}) by mathematical induction:

\smallskip

\noindent
\underline{Case $m=0$}. \quad
Since $|u_-| \le \sqrt{\frac{|p'(v_+)|}{2}}v_+$, which easily follows from (\ref{econdi1}),
we have 
\begin{align*}
\label{m0} 
|r^{2(n-1)}{\eta^{E}}^{(0)}(r)|
\le r^{2(n-1)}\int_r^{\infty}\frac{(n-1)|u_-|^2 }{\frac{1}{2}|\tilde p'(v_+)|s^{2n-1}v_+}ds
\le \frac{|u_-|^2}{|p'(v_+)|v_+},\quad r \ge 1.
\end{align*}
Hence, we ask the constant 
$M$ to satisfy 
\begin{align}
\frac{|u_-|^2}{|p'(v_+)|v_+} \le M,
\end{align}
so that it holds $\|{\eta^{E}}^{(0)}\|_X \le M$. 

\smallskip

\noindent
\underline{Case $m=k+1$}\ $(k\ge 0)$.\quad
Suppose $\|{\eta^{E}}^{(k)}\|_X \le M$. Here we also ask the constant $M$ to satisfy another assumption
\begin{align}
M\le \frac{v_+}{2},
\end{align}
which, in particular, implies
\begin{align}
\label{bound}
|{\eta^{E}}^{(k)}(r)|\le \frac{v_+}{2},\quad 
|{\e^{E}}^{(k)}(r)|=|\frac{u_-}{v_++{\eta^{E}}^{(k)}(r)}|\le |\frac{2u_-}{v_+}|,
\quad r \ge 1.
\end{align}
Then, we estimate ${\eta^{E}}^{(k+1)}$ defined by (\ref{eas}) as follows:
\begin{align}
\begin{aligned}
&|r^{2(n-1)}{\eta^{E}}^{(k+1)}(r)|\\
&\le |r^{2(n-1)} \int_r^{\infty}\frac{(n-1)(v_++{\eta^{E}}^{(k)}(s))|\frac{2u_-}{v_+}|^2}
{\left(|\tilde p'(\frac{3}{2}v_+)|-|\frac{2u_-}{v_+}|^2 \right)s^{2n-1}}ds|\\ 
& \le |r^{2(n-1)}\frac{2(n-1)|\frac{2u_-}{v_+}|^2}{|\tilde p'(\frac{3}{2}v_+)|} 
\int_r^{\infty}\frac{(v_++{\eta^{E}}^{(k)}(s))}{s^{2n-1}}ds|  
\le  \frac{|\frac{2u_-}{v_+}|^2(v_++M)}{|\tilde p'(\frac{3}{2}v_+)|},
\end{aligned} 
\end{align}
where we used the assumption (\ref{econdi0}) and (\ref{bound}). 
Therefore we obtain 
\begin{align}
\label{ek+1e}
\|{\eta^{E}}^{(k+1)}\|_X
\le \frac{|\frac{2u_-}{v_+}|^2(v_++M)}{|\tilde p'(\frac{3}{2}v_+)|}.
\end{align}
Here, we further assume 
\begin{equation}
\label{econdi2}
\frac{4|u_-|^2v_+}{|v_+|^2|\tilde p'(\frac{3}{2}v_+)|}\le \frac{M}{2},\qquad 
\frac{4|u_-|^2}{|v_+|^2|\tilde p'(\frac{3}{2}v_+)|}\le \frac{1}{2},
\end{equation}
so that  (\ref{ek+1e}) gives the desired estimate $\|{\eta^{E}}^{(k+1)}\|_X \le M$. By elementary calculations, it is easy to see that there exists a
positive constant $\epsilon_0$ such that if 
$|u_-| \le \epsilon_0$, the assumptions 
(\ref{econdi1}), (3.12),(3.13) and (\ref{econdi2}) hold, and in particular,
$M$ can be chosen by
\begin{equation*}
M = \frac{8|u_-|^2v_+}{|v_+|^2|\tilde p'(\frac{3}{2}v_+)|},
\end{equation*}
which proves the uniform boundedness of 
${\eta^{E}}^{(m)}\ (m\ge 0)$ in $X$ with the estimate (\ref{eu}).
Once (\ref{eu}) is proved, the proof to show $\{{\eta^{E}}^{(m)}\}_{m\ge 0}$ is a Cauchy sequence in $X$ is very standard. In fact, we may estimate 
\begin{align*}
\begin{aligned}
\|{\eta^{E}}^{(m+1)}-{\eta^{E}}^{(m)}\|_X
&=\sup_{r\ge 1}|r^{2(n-1)}\int_{r}^{\infty}(G[{\eta^{E}}^{(m)}]-G[{\eta^{E}}^{(m-1)}])(s)\,ds|\\[2mm]
&\le C|u_-|^2 \|{\eta^{E}}^{(m)}-{\eta^{E}}^{(m-1)}\|_X, \quad m \ge 1,
\end{aligned}
\end{align*}
in the same way as in (\ref{m0})-(\ref{ek+1e}), and taking $|u_-|$ suitably small again if needed, 
we can show
\begin{align*}
\|{\eta^{E}}^{(m+1)}-{\eta^{E}}^{(m)}\|_X
\le \frac{1}{2}\|{\eta^{E}}^{(m)}-{\eta^{E}}^{(m-1)}\|_X, \quad m \ge 1,
\end{align*}
which proves that $\{{\eta^{E}}^{(m)}\}_{m\ge 0}$
is a Cauchy sequence in the function space $X$. 
Thus, as the limit for $m$ in the approximate sequence (\ref{eas}), 
the solution $\eta^E$ of (\ref{E2}) with the desired estimate (\ref{decayE1}) is obtained. 
The arguments on the regularity and uniqueness of the solution 
are also very standard, so we omit them. 
The decay rate of $\eta^{E}_r$ in (\ref{decayE2}) can also be directly derived 
by the formula (\ref{E1}). \quad $\Box$

\subsection{Zero viscosity limit for the outflow problem }
In this subsection, we consider 
the zero viscosity limit for the outflow problem.
%to the compressible Navier-Stokes Equation. 
First, we rewrite the equation (\ref{reformE}) as in the form
\begin{align}
\label{out1}
 \left\{
\begin{aligned}
&-\frac{\big(\tilde p(v_++\eta^{E})- \tilde p(v_+)\big)}{\epsilon^{E} \mu}r^{n-1}
-\frac{{\epsilon^{E}} v_+}{2\mu}\frac{1}{r^{n-1}}\\[2mm]
&\qquad\qquad\qquad
-\frac{{\epsilon^{E}}\eta^{E}}{\mu r^{n-1}}+\frac{{\epsilon^{E}}(n-1)r^{n-1}}{\mu}\int_r^{\infty}\frac{\eta^{E}(s)}{s^{2n-1}}ds=0,\quad r \ge 1,\\[2mm]
&\displaystyle{\lim_{r\to \infty}}\eta^{E}(r) =0.
\end{aligned}
 \right.
\end{align}
\noindent
Now we introduce a new variable $\chi(r)$ as $\chi(r):=\eta(r)-\eta^{E}(r)$. 
Subtracting (\ref{out1}) from (\ref{reform}) 
and noting the definitions $\epsilon = u_-/(v_++\eta(1))$ and $\epsilon^{E} = u_-/(v_++\eta^{E}(1))$, 
%defined in (\ref{Nep}) and (\ref{Eep}), respectively, 
we obtain the equation of $\chi(r)$ as 
\begin{align}
\label{out2}
 \left\{
\begin{aligned}
&\chi_r-\frac{r^{n-1}}{\epsilon\mu}a(r)\chi= -\eta^{E}_r+F[\chi](r),\quad r \ge 1,\\[2mm]
&\displaystyle{\lim_{r\to \infty}}\chi(r) =0,
\end{aligned}
 \right.
\end{align}
where
$$
a := \frac{\tilde p(v_++\eta^{E}+\chi)-\tilde p(v_++\eta^{E})}{\chi},
$$
and
\begin{align*}
\begin{aligned}
%a(\eta) &:= \frac{\tilde p(v_++\eta^{E}+\chi)-\tilde p(v_++\eta^{E})}{\chi},\\
%and, & \\
F[\chi](r) &:=\frac{r^{n-1}\chi(1)}{\mu u_-}\left(\tilde p(v_++\eta^{E})-\tilde p(v_+)\right)
-\frac{\chi(1)u_-}{(\eta^{E}(1)+v_+)(\eta(1)+v_+)}\frac{v_+}{2\mu r^{n-1}}\\[2mm]
&+\frac{u_-\chi(r)}{\mu r^{n-1}(\eta(1)+v_+)}-\frac{\chi(1)u_-}{(\eta^{E}(1)+v_+)(\eta(1)+v_+)}
\frac{\eta^{E}(r)}{\mu r^{n-1}}\\[2mm]
&-\frac{u_-(n-1)r^{n-1}}{\mu(\eta(1)+v_+)}\int_r^{\infty}\frac{\chi(s)}{s^{2n-1}}ds
+\frac{\chi(1)u_-(n-1)r^{n-1}}{\mu (\eta^{E}(1)+v_+)(\eta(1)+v_+)}
\int_r^{\infty}\frac{\eta^{E}(s)}{s^{2n-1}}ds\\[2mm]
&=: I_1+I_2+I_3+I_4+I_5+I_6.
\end{aligned}
\end{align*}
Here we note that we already know the existence of the solution 
$\eta(r)$ of (\ref{reform}) and $\eta^{E}(r)$ of (\ref{E2}) with the estimates (\ref{estouteta}) and (\ref{decayE1}) respectively. 
Applying  the Duhamel's principle to the equation (\ref{out2}) in terms of $\chi$, we obtain 
\begin{align}
\label{intchi}
\begin{aligned}
\chi(r)&=-\int_r^{\infty}e^{-\int_{r}^{s}\frac{a(\tau)}{\epsilon\mu}\tau^{n-1}d\tau}(-\eta^{E}_r+F[\chi])(s)\,ds\\
&=-\int_r^{\infty}e^{-\int_{r}^{s}\frac{a(\tau)}{\epsilon\mu}\tau^{n-1}d\tau}
(-\eta^{E}_r+I_1+I_2+I_3+I_4+I_5+I_6)(s)\,ds.
\end{aligned}
\end{align}
By the assumption $\tilde p'(v)<0\ (v>0)$, Theorem 2.1, and (\ref{decayE1}), 
it is easy to see that if $|u_-|+|v_--v_+|$ is suitably small, there exist a positive constant 
$\delta$ 
which is independent of $\mu$
satisfying
\begin{equation*}
-a(r) \ge \delta, \qquad r\ge 1.
\end{equation*}
Therefore by using the estimates for $\eta$ and $\eta^{E}$ already obtained, we can estimate the first and second term on the right hand side of (\ref{intchi}) as 
\begin{align}
\label{chiu1}
\begin{aligned}
|\int_r^{\infty}e^{-\int_{r}^{s}\frac{a(\tau)}{\epsilon\mu}\tau^{n-1}d\tau}&\eta^{E}_r(s)\, ds| \\
\le  &\int_r^{\infty}e^{-\frac{v(1)\delta}{|u_-|\mu}(s-r)}|\eta^{E}_r(s)|\, ds \\
\le &\int_r^{\infty}e^{-\frac{v(1)\delta}{|u_-|\mu}(s-r)}ds 
\cdot \sup_{r\ge 1}|\eta^{E}_r(r)| \\
\le &\frac{|u_-|\mu}{v(1)\delta}\cdot \sup_{r\ge 1}|\eta^{E}_r(r)|
\le C \mu; \\
|\int_r^{\infty}e^{-\int_{r}^{s}\frac{a(\tau)}{\epsilon\mu}\tau^{n-1}d\tau}&I_1(s)\,ds| \\
\le  &\frac{|\chi(1)|}{\mu |u_-|}\int_r^{\infty}e^{-\frac{v(1)\delta}{|u_-|\mu}(s-r)}|s^{n-1}
\left(\tilde p(v_++\eta^{E})-\tilde p(v_+)\right)| \,ds \\
\le  &\frac{C|\chi(1)|}{\mu |u_-|}\int_r^{\infty}e^{-\frac{v(1)\delta}{|u_-|\mu}(s-r)}ds \cdot \sup_{r\ge 1}|r^{n-1} \eta^{E}(r)| \\
\le &\frac{C|u_-|^2}{v(1)\delta}\cdot \sup_{r\ge 1}|\chi(r)|
\le C|u_-|\sup_{r\ge 1}|\chi(r)|.
\end{aligned}
\end{align}
By applying the same manner to the remaining terms on the right hand side of (\ref{intchi}), we can obtain
%\begin{align}
%\label{chiu2}
%\begin{aligned}
%&\int_r^{\infty}e^{-\int_{r}^{s}\frac{a(\tau)}{\epsilon\mu}\tau^{n-1}d\tau}I_1(s) ds \\
%&\le C\frac{\chi(1)}{\mu u_-}\int_r^{\infty}e^{-\int_{r}^{s}\frac{a(\tau)}{\epsilon\mu}\tau^{n-1}d\tau}s^{n-1}\eta^{E}(s) ds \le \frac{C|u_-|^2}{v(1)\delta}\sup_{r\ge 1}|\chi(r)|
%\end{aligned}
%\end{align}
\begin{align}
\label{chiu3}
\begin{aligned}
&|\int_r^{\infty}e^{-\int_{r}^{s}\frac{a(\tau)}{\epsilon\mu}\tau^{n-1}d\tau}
(I_2+I_3+I_4+I_5+I_6)(s)\, ds| \le C|u_-|\sup_{r\ge 1}|\chi(r)|.
\end{aligned}
\end{align}
Substituting (\ref{chiu1}) and (\ref{chiu3}) into (\ref{intchi}), we can estimate $|\chi(r)|$ as 
\begin{align}
\label{chiu4}
\sup_{r\ge 1}|\chi(r)|\le C\mu+C|u_-|\sup_{r\ge 1}|\chi(r)|.
\end{align}
Hence, taking $|u_-|$ suitably small on the right hand side of (\ref{chiu4}), 
we obtain
\begin{align}
\label{chiu5}
\sup_{r\ge 1}|\chi(r)|=
\sup_{r\ge 1}|(\eta-\eta^{E})(r)|
\le C\mu. 
\end{align}
Finally, recalling $u$ and $u^{E}$ are respectively given by
$$
u(r)=\frac{u_-(v_++\eta(r))}{r^{n-1}(v_++\eta(1))},\quad
u^{E}(r)=\frac{u_-(v_++\eta^{E}(r))}{r^{n-1}(v_++\eta^{E}(1))},\quad r\ge 1,
$$
we can also show
$$
\sup_{r\ge 1}|(u-u^{E})(r)| \le C\mu. 
$$
This completes the proof for the result (II) in Theorem \ref{etE}. \qquad $\Box$

\section{Inflow problem}
In this section, we consider the case $u_->0$ and show the statement (III)
in Theorem \ref{bl+E}. That is, we show that the inviscid limit of the solution of 
the Navier-Stokes equation (\ref{nrs0}) is given by
a linear  superposition of that of the Euler equation (\ref{Es}) 
with boundary condition $\rho^{E}(r_0)u^{E}(r_0)=\rho_-u_-$ 
and the boundary layer solution of (\ref{Ebl}).
%%%%%%%%%%%%%%%%%%%%%%%%%%%%%%%%%%%%%%%%%%%%%%%%%%%%%%%%

%%%%%%%%%%%%%%%%%%%%%%%%%%%%%%%%%%%%%%%%%%%%%%%%%%%%%%%%
\subsection{Euler equation} 
In this subsection, we only make a short comment on the existence of 
solution of the Euler equation (\ref{Es}) which corresponds to the inflow problem,
since the arguments are almost the same as in the subsection 3.1.
For the case of the inflow problem,
we naturally look for the solution $(\rho^{E},u^{E})$ of (\ref{Es}) 
satisfying 
\begin{align}
\label{4n1e}
r^{n-1}\rho^{E}(r) u^{E}(r)=\epsilon := \rho_- u_-, \qquad r\ge 1.
\end{align}
%where we note that $\epsilon$ in (\ref{4n1e}) is 
%the constant defined in (\ref{n1}), that is $\epsilon=\rho_- u_-$. 
To construct the solution, we again introduce a new variable 
$\eta^{E}(r):=v^{E}(r)-v_+$. Then, due to the same argument as in 
deriving (\ref{uE})-(\ref{E1}),
we may consider the solution of the following integral equation 
%On the other hand, by the equation in (\ref{4n1e}), the velocity $u^{E}$ is given in terms of specific volume $v^{E}(r):=1/\rho^{E}(r)$ as
%\begin{align}
%\label{4uE}
%u^{E}(r)=\frac{\epsilon}{r^{n-1}}v^{E}(r),\quad r\ge 1.
%\end{align}
%Substituting (\ref{uE}) into the second equation of (\ref{Es}) and introducing a new variable $\eta^{E}(r):=v^{E}(r)-v_+$ we have
%\begin{align}
%\label{4baseE}
%\frac{{\epsilon}^2}{r^{n-1}}\left(\frac{v_++\eta^{E}}{r^{n-1}}\right)_r+\tilde{p}(v_++\eta^{E})_r=0, 
%\end{align}
%which also can be rewritten as
%\begin{align}
%\label{4E1}
%\eta^{E}_r(r)=\frac{{\epsilon}^2(n-1)(v_++\eta^{E})}{(|\epsilon^{E}|^2r^{-2(n-1)}+\tilde p'(v_++\eta^{E}))r^{2n-1}}.
%\end{align}
%Integrating (\ref{E1}) in terms of space variable on $[r,\infty]$, we obtain 
\begin{align}
\label{4E2}
\begin{aligned}
\eta^{E}(r)&=\int_r^{\infty}\frac{{\epsilon}^2(n-1)(v_++\eta^{E})}
{(|\tilde p'(v_++\eta^{E})|-{\epsilon}^2s^{-2(n-1)})s^{2n-1}}\,ds.\\
%&:=\int_r^{\infty}G[\eta^{E}](r)ds,
\end{aligned}
\end{align}
%\begin{align}
%\label{E2}
%\eta^{E}(r)=-\int_r^{\infty}\frac{|\epsilon^{E}|^2(n-1)(v_++\eta^{E})}{(|\epsilon^{E}|^2s^{-2(n-1)}+\tilde p'(v_++\eta^{E}))s^{2n-1}}ds,
%\end{align}
Moreover, adopting the same strategy as (\ref{base})-(\ref{reform}), 
we see that (\ref{4E2}) is equivalent to 
%\begin{align}
%\label{4reformE}
% \left\{
%\begin{aligned}
%&\big(\tilde p(v_++\eta^{E})- \tilde p(v_+)\big)
%+\frac{{\epsilon}^2 v_+}{2}\frac{1}{r^{2(n-1)}}
%+\frac{{\epsilon}^2\eta^{E}}{r^{2(n-1)}}-{\epsilon}^2(n-1)\int_r^{\infty}\frac{\eta^{E}(s)}{s^{2n-1}}ds=0,\\[2mm]
%&\displaystyle{\lim_{r\to \infty}}\eta^{E}(r) =0,\quad r>1. 
%\end{aligned}
% \right.
%\end{align}
\begin{align}
\label{4reformE}
 \left\{
\begin{aligned}
&\big(\tilde p(v_++\eta^{E})- \tilde p(v_+)\big)
+\frac{{\epsilon}^2 v_+}{2}\frac{1}{r^{2(n-1)}}\\[2mm]
&\qquad\qquad\qquad\qquad\quad
+\frac{{\epsilon}^2\eta^{E}}{r^{2(n-1)}}-{\epsilon}^2(n-1)\int_r^{\infty}\frac{\eta^{E}(s)}{s^{2n-1}}ds=0,
\quad r \ge 1,
\\[2mm]
&\displaystyle{\lim_{r\to \infty}}\eta^{E}(r) =0.
\end{aligned}
 \right.
\end{align}
By the same computation in the proof of Theorem \ref{rmtE}, 
we can obtain the following theorem for the reformulated problem (\ref{4E2})/(\ref{4reformE}).
\begin{thm} 
\label{4rmtE}
Let $n \ge 2$. Then, for any $v_+>0$, there exist positive constants $\epsilon_0$ and $C$ 
such that if $|u_-|+|\rho_--\rho_+| \le \epsilon_0$, 
there exists a unique smooth solution $\eta^{E}$ 
of the problem {\rm (\ref{4E2})/(\ref{4reformE})} satisfying
\begin{subequations}
\begin{align}
\label{4decayE1}
&|\eta^{E}(r)|\le Cr^{-2(n-1)}|u_-|^2, \qquad r\ge 1, \\[2mm]
\label{4decayE2}
&|\eta^{E}_r(r)|\le Cr^{-2n+1}|u_-|^2, \qquad \ \ r\ge 1.
\end{align}
\end{subequations}
\end{thm}
Once Theorem 4.1 is proved, since the velocity $u^{E}$ is given by 
the relation (4.1),
thus the statement (I) in Theorem 1.3 is also proved.
%\noindent
%The arguments of the proof of theorem \ref{4rmtE} is 
%almost the same as that of Theorem \ref{rmtE}, so we omit that.

\subsection{Boundary layer solution} 

In this subsection, we construct a boundary layer equation 
and its solution. 
We derive this by introducing a new independent 
variable $y$ and new function $\tilde \eta(y)$  as 
\begin{align}
\label{bltrans}
y:=\frac{r-1}{\mu}, \qquad \tilde \eta(y):=\eta(1+\mu y).
\end{align}
By applying (\ref{bltrans}) to the arranged Navier-Stokes equation (\ref{Narrange}), 
we can rewrite (\ref{Narrange}) in terms of $y$ and $\tilde \eta(y)$ as
\begin{align}
\label{Narrange1}
\begin{aligned}
&\epsilon \left(\frac{\tilde{\eta}_y}{(1+\mu y)^{n-1}}\right)_y=\tilde p(v_++\tilde \eta)_y \\
&\quad +\frac{\epsilon^2 v_+}{2}\Big(\frac{1}{(1+\mu y)^{2(n-1)}}\Big)_y
+\frac{\epsilon^2}{(1+\mu y)^{n-1}}\left(\frac{\tilde \eta}{(1+\mu y)^{n-1}}\right)_y,\qquad y\ge 0.
\end{aligned}
\end{align}
Let $\mu=0$ in (\ref{Narrange1}) and define $\eta^B$ 
as a solution of limit of the equation (\ref{Narrange1}). This computation gives 
\begin{align}
\label{Narrange2}
\epsilon \eta^B_{yy}=\tilde p(v_++\eta^B)_y+\epsilon^2 \eta^B_y.
\end{align}
Integrating (\ref{Narrange2}) with respect to $y$ on $[y, \infty)$ 
under the assumption $\displaystyle{\lim_{y\to \infty}}\eta^B(y) =\eta^{E}(1)$, 
which is the boundary data of the solution  of 
the Euler equation obtained in the subsection 3, 
we obtain 
%the problem for boundary layer equation as 
\begin{align}
\label{Narrange3}
 \left\{
\begin{aligned}
  &\eta^B_y= 
\frac{1}{\epsilon}\big(\tilde p(v_++\eta^B(y))- \tilde p(v_++\eta^{E}(1))\big)
+\e (\eta^B(y)-\eta^{E}(1)), \quad {y \ge 0}, \\[2mm]
 &\displaystyle{\lim_{y\to \infty}}\eta^B(y) =\eta^{E}(1), 
 \quad \eta^B(0)= \eta_-:=v_--v_+,
 \end{aligned}
 \right.
\end{align}
where $\e=u_-/v_-$. 
%In this paper, we define the ``boundary layer solution" 
%as the solution of the problem (\ref{Narrange3}).
In this paper, we call the equation (\ref{Narrange3}) 
``boundary layer equation", and its solution 
``boundary layer solution".

Next we show the existence of the boundary layer solution of (\ref{Narrange3}). 
The theorem for the problem (\ref{Narrange3}) which we need to prove is
\begin{thm} 
\label{blexist}
Let $n \ge 2$. Then, for any $v_+>0$, there exist positive constants 
$\epsilon_0$, $\nu_0$ and $C$ such that if $|u_-|+|\eta_--\eta^{E}(1)| \le \epsilon_0$, 
there exists a unique smooth solution $\eta^B$ of the problem {\rm (\ref{Narrange3})} satisfying
\begin{align}
\label{4decaybl}
|\eta^B(y)-\eta^{E}(1)|\le C|\eta_--\eta^{E}(1)|
\ e^{-\frac{\nu_0}{\e}y}, \qquad y\ge 0. 
\end{align}
\end{thm}
Once Theorem 4.2 is proved, noting that 
$$
{\rho}^B(y) := \frac{1}{v_++\eta^B(y)},\quad
{u}^B(y) := \frac{u_-}{v_-}(v_++\eta^B(y))
$$
give the solution of (1.11), and that Theorem 4.1 gives
$|\eta^{E}(1)| \le C|u_-|^2$, 
we can easily show the statement (II) in Theorem 1.3.

\medskip

{\it Proof of Theorem \ref{blexist}}\quad  We first rewrite (\ref{Narrange3}) as 
\begin{align}
\label{Narrange4}
\begin{aligned}
  &(\eta^B-\eta^{E}(1))_y-\frac{\tilde p'(v_*)}{\e}(\eta^B-\eta^{E}(1))\\
&=\frac{1}{\epsilon}\left\{\tilde p(v_++\eta^B(y))- \tilde p(v_*)-\tilde p'(v_*)(\eta^B-\eta^{E}(1))\right\}+\e (\eta^B(y)-\eta^{E}(1)),
 \end{aligned}
\end{align}
where $v_*=v_++\eta^{E}(1)$.
Now, we introduce a new unknown function $\bar \eta$, as the deviation of $\eta^B$ from the far field state $\eta^{E}(1)$ by
\begin{equation}
\label{bareta}
\hat \eta(y) := \eta^B(y) - \eta^{E}(1),\qquad y \ge 0.
\end{equation}
Plugging (\ref{bareta}) into (\ref{Narrange4}), we obtain 
\begin{align}
\label{Narrange5}
 \left\{
\begin{aligned}
 &\hat \eta_y+a_\e \hat \eta=\frac{1}{\e}N[\hat{\eta}], \quad y\ge 0,\\[2mm]
 &\displaystyle{\lim_{y\to \infty}}\hat \eta(y) =0, \quad \hat \eta(0)= \eta_- - \eta^{E}(1),
 \end{aligned}
 \right.
\end{align}
where 
\begin{align*}
\begin{aligned}
&a_\e := -\frac{\tilde p'(v_*)}{\e}-\e,\\[2mm]
&N[\hat{\eta}]:=\tilde p(v_*+\hat{\eta})- \tilde p(v_*)-\tilde p'(v_*)\hat{\eta}.
\end{aligned}
\end{align*}
Applying the Duhamel's principle to the equation
(\ref{Narrange5}) in terms of $\hat \eta$, we obtain
\begin{align}
\label{Narrange6}
\hat \eta(y)=e^{-a_{\e}y}\hat \eta(0)
+\frac{1}{\e}\int_0^y e^{-a_{\e}(y-s)}N[\hat \eta](s)\,ds.
\end{align}
\noindent
Here, by $\tilde{p}'(v_+)<0$ and $|\eta^{E}(1)|\le C|u_-|^2$, it easy to see that
for any $v_+>0$ there exist the positive constants $\e_0$ and $\nu_0$
such that if $|v_--v_+| +|u_-|\le \e_0$, it holds 
$a_\e \ge \nu_0/{\e}$.
Then, we may look for the solution of (\ref{Narrange6}) in a
Banach space $Y$, with its norm $\|\cdot\|_Y$, defined by
\begin{align*}
%\label{Y}
Y:=\{g \in C([1,\infty));\ \sup_{y\ge 0}|e^{\frac{\nu_0}{\e}y} g(y)|< \infty \},
\quad \|g\|_Y = \sup_{y\ge 0}|e^{\frac{\nu_0}{\e}y}g(y)|.
\end{align*}
Since  $N[\hat \eta]= O({\hat \eta}^2)$,
the remaining arguments to obtain the desired solution of (\ref{Narrange6}) in $Y$
for suitably small $\hat{\eta}(0)$ are quite standard and much simpler than
that in the proof of Theorem 3.1.
So we omit the details. $\Box$

\subsection{Zero viscosity limit for the inflow problem }
In this subsection, we show the statement (III) in Theorem \ref{bl+E}.
We first note that if we change the variable $y$ back to $r$ for the
boundary layer solution ${\eta}^B(y)$, and introduce $\bar \eta(r)$ by
$$
{\bar \eta}(r) = {\eta}^B(\frac{r-1}{\mu}) - \eta^{E}(1),
$$
the boundary layer equation  (\ref{Narrange3}) is rewritten in terms of ${\bar \eta}$
as in the form
\begin{align}
\label{bl_r}
\bar{\eta}_r-\frac{\tilde p(v_++\eta^{E}(1)+\bar\eta)
-\tilde p(v_++\eta^{E}(1))}{\e\mu}-\frac{\e}{\mu}\bar\eta=0,
\end{align}
where recall $\epsilon = u_-/v_->0$.  
Then, to show the statement (III) in Theorem \ref{bl+E} on $\rho$,  
we may show the solution $\eta$ of (\ref{reform})
uniformly converges to the linear combination of 
the solutions of (\ref{4reformE}) and (\ref{bl_r}), that is,
$\eta^{E}+\bar \eta$, as
$\mu \to 0$. To do that,
keeping in mind that for any positive constant $\al\in (0,1)$
it holds $|\bar \eta(r)|\le C\exp\{-\nu_0\mu^{\al-1}/\e\}$ for $r\ge 1+\mu^\al$ by Theorem 4.2,
we divide the interval $[1,\infty)$ into $I_b:=[1,1+\mu^{\al}]$ 
and $I_{\infty}:=[1+\mu^{\al},\infty)$. Then,
we first show that on the interval $I_{\infty}$,
 $\eta$ uniformly converges toward $\eta^{E}$, and next
 on the interval $I_{b}$, toward $\eta^{E}+\bar\eta$,
 as $\mu \to 0$.

%the solution of (\ref{reform}) tends to the solution of Euler equation of (\ref{Es}) 
%as $\mu\to 0$. On the other hand we show that 
%on interval $I_{b}$, the solution tends to the superposition 
%of boundary layer solution of (\ref{Ebl}) and the solution 
%of Euler equation (\ref{Es}) as $\mu\to 0$.

\medskip

\noindent
\underline{\bf zero viscosity limit on $I_{\infty}$}

\smallskip

We define a new variable $\chi(r)$ as the difference of the solution $\eta$ of
 Navier-Stokes equation  (\ref{reform}) from the solution $\eta^{E}$ 
 of Euler equation  (\ref{4reformE}), that is  $\chi(r):=\eta(r)-\eta^{E}(r)$. 
To have the equation of $\chi(r)$, we rewrite the equation (\ref{4reformE}) as in the form
\begin{align}
\label{ineuler}
 \left\{
\begin{aligned}
&-\frac{\big(\tilde p(v_++\eta^{E})- \tilde p(v_+)\big)}{\epsilon \mu}r^{n-1}
-\frac{{\epsilon} v_+}{2\mu}\frac{1}{r^{n-1}}\\[2mm]
&\qquad\qquad\qquad
-\frac{{\epsilon}\eta^{E}}{\mu r^{n-1}}
+\frac{{\epsilon}(n-1)r^{n-1}}{\mu}\int_r^{\infty}\frac{\eta^{E}(s)}{s^{2n-1}}ds=0,\\[2mm]
&\displaystyle{\lim_{r\to \infty}}\eta^{E}(r) =0.
\end{aligned}
 \right.
\end{align}
\noindent
Subtracting (\ref{ineuler}) from (\ref{reform}), we obtain the equation of $\chi(r)$ as 
\begin{align}
\label{inf0}
 \left\{
\begin{aligned}
&\chi_r+\frac{r^{n-1}}{\epsilon\mu}a(r)\chi= 
-\eta^{E}_r-\frac{\e (n-1)r^{n-1}}{\mu}\int_r^{\infty}\frac{\chi(s)}{s^{2n-1}}ds,\quad r > 1,\\[2mm]
&\displaystyle{\lim_{r\to \infty}}\chi(r) =0,\quad \chi(1) = \eta_--\eta^{E}(1),
\end{aligned}
 \right.
\end{align}
where
\begin{align*}
\begin{aligned}
a := -\frac{\tilde p(v_++\eta^{E}+\chi)-\tilde p(v_++\eta^{E})}{\chi}-\frac{\e^2}{r^{2(n-1)}}.
\end{aligned}
\end{align*}
By the assumption $\tilde p'(v)<0\ (v>0)$, Theorem 2.1, and (\ref{4decayE1}), 
it is easy to see that there exist a positive constant $\delta$ which is independent of
$\mu$ satisfying
\begin{equation}
\label{|a|}
a(r) \ge \delta, \qquad r\ge 1,
\end{equation}
for suitably small $|v_--v_+|$ and $|u_-|$.
The theorem for the reformulated problem (\ref{inf0}) is as follows.
\begin{thm} 
\label{ivli0}
Let $n \ge 2$ and $\al \in (0,1)$. Then, for any $v_+>0$, there exists positive constants $\epsilon_0$ 
and $C$ which are independent of $\mu\in (0,1]$
such that if $|u_-|+|v_--v_+| \le \epsilon_0$, 
the solution $\chi$  of (\ref{inf0}) satisfies
\begin{align}
\label{invsidlimitinf}
\sup_{r\ge 1+\mu^{\alpha}}|\chi(r)| \le C\mu.
\end{align}
\end{thm}

{\it Proof.} \ Applying the Duhamel's principle to the equation (\ref{inf0}) 
in terms of $\chi$, we obtain 
\begin{align}
\label{intchiinf}
\begin{aligned}
\chi(r)&=\chi(1)e^{-\int_1^r \frac{a(\tau)}{\e\mu}\tau^{n-1}d\tau}
-\int_1^r e^{-\int_{\tau}^{r}\frac{a(s)}{\epsilon\mu}s^{n-1}ds}\eta^{E}_r({\tau})\,d\tau
\\
&\quad\quad -\frac{\e (n-1)}{\mu}\int_1^r e^{-\int_{\tau}^{r}\frac{|a(s)|}{\epsilon\mu}s^{n-1}ds}
(\tau^{n-1}\int_{\tau}^{\infty}\frac{\chi(s)}{s^{2n-1}}ds)\,d\tau \\[2mm]
&\quad\quad =: I_1+I_2+I_3.
\end{aligned}
\end{align}
By using (\ref{|a|}) and noting $r\ge 1+\mu^{\alpha}$, we can estimate $I_1$ as 
\begin{align}
\label{I1bleu}
\begin{aligned}
|I_1| &\le |\chi(1)|e^{-\int_1^r \frac{\delta}{\e\mu}\tau^{n-1} d\tau}
\le |\chi(1)|e^{-\frac{\delta}{\e \mu}(r-1)}\\
%\le |\chi(1)|e^{-\frac{\delta}{\e} \mu^{\al -1}}\\
&\le |\chi(1)|e^{-\frac{\delta}{\e} \frac{1}{\mu^{1-\al}}}
\le C\mu.
\end{aligned}
\end{align}
By Theorem 2.1 and (\ref{|a|}), we can estimate $I_2$ in (\ref{intchiinf}) as 
\begin{align}
\label{I2bleu}
\begin{aligned}
|I_2| &\le \sup_{r\ge 1}|\eta^{E}_r(r)|\int_1^r e^{-\int_{\tau}^{r}\frac{\delta}{\epsilon\mu}s^{n-1}ds}d\tau 
\le \frac{\e\mu}{\delta}\sup_{r\ge 1}|\eta^{E}_r(r)| \le C\mu.
\end{aligned}
\end{align}
By (\ref{|a|}) and integration by parts, we can estimate $I_3$  in (\ref{intchiinf}) as 
\begin{align}
\label{I3bleu}
\begin{aligned}
|I_3|
&\le \frac{\e (n-1)}{\mu}\int_1^r e^{-\frac{\delta}{\e\mu n}(r^n-\tau^n)}
\tau^{n-1}\int_{\tau}^{\infty}\frac{|\chi(s)|}{s^{2n-1}}ds\,d\tau \\
&= \frac{\e^2 (n-1)}{\delta}\int_1^r \left(e^{-\frac{\delta}{\e\mu n}(r^n-\tau^n)}\right)_{\tau}\int_{\tau}^{\infty}\frac{|\chi(s)|}{s^{2n-1}}dsd\tau \\
&\le \frac{\e^2 (n-1)}{\delta}\left(\int_{r}^{\infty}\frac{|\chi(s)|}{s^{2n-1}}ds+\int_1^r  e^{-\frac{\delta}{\e\mu n}(r^n-\tau^n)} \frac{|\chi(\tau)|}{\tau^{2n-1}}d\tau \right) \\
&\le \frac{\e^2 (n-1)}{\delta}\left(C \sup_{r\ge 1+\mu^{\alpha}}|\chi(r)|+C\frac{\e\mu}{\delta}\sup_{r\ge 1}|\chi(r)|\right)\\
&\le C|u_-| \sup_{r\ge 1+\mu^{\alpha}}|\chi(r)|+C\mu.
\end{aligned}
\end{align}
Substituting (\ref{I1bleu})-(\ref{I3bleu}) into (\ref{intchiinf}) 
and taking $u_-$ suitably small, 
we finally obtain
$$
\sup_{r \ge 1+\mu^{\alpha}}|\chi(r)| \le C\mu.
$$
Thus the proof of Theorem 4.3 is completed.
%\begin{align}
%\label{4supmigi}
%\begin{aligned}
%\sup_{1+\mu^{\alpha}\le r}|\chi(r)| 
%&\le C\left(|\chi(1)|e^{-\frac{\delta}{\e} \frac{1}{\mu^{1-\al}}}
%+\frac{\e\mu}{\delta}\sup_{1\le r}|\eta^{E}_r|+\frac{\mu\e^3 (n-1)}{\delta^2}\sup_{1\le r}|\chi(r)|\right).
%%& \to 0 \quad (\mu \to 0).
%\end{aligned}
%\end{align}
%Finally taking limit of $\mu$ to $0$, the right hand side of (\ref{4supmigi}) tends to $0$, that is 
%\begin{align}
%\label{limitchir}
%\begin{aligned}
%\sup_{1+\mu^{\alpha}\le r}|\chi(r)| \to 0, \quad (\mu \to 0),
%\end{aligned}
%\end{align}
$\Box$

\medskip

\noindent
\underline{\bf zero viscosity limit on $I_{b}$}

\smallskip

We define a new variable $\chi(r)$ as the difference of the solution $\eta$ 
of the Navier-Stokes equation  (\ref{reform}) from the linear superposition 
of $\eta^{E}$ and $\bar \eta$, which are the solutions of the Euler equation  (\ref{4reformE}) 
and boundary layer equation (\ref{bl_r}), 
that is, $\chi(r):=\eta(r)-(\eta^{E}(r)+\bar \eta(r))$. 
To have the equation in terms of $\chi(r)$, 
%\begin{align}
%\label{bl_r}
%\bar{\eta}_r-\frac{\tilde p(v_++\eta^{E}(1)+\bar\eta)-\tilde p(v_++\eta^{E}(1))}{\e\mu}-\frac{\e}{\mu}\bar\eta=0.
%\end{align}
subtracting (\ref{ineuler}) and (\ref{bl_r}) from (\ref{reform}) gives
%\begin{align}
%\begin{aligned}
%&\chi_r+\eta^{E}_r-\frac{r^{n-1}}{\e\mu}(\ti p(v_++\eta^{E}+\bar \eta+\chi)-\ti p(v_++\eta^{E}+\bar \eta))\\
%&-\frac{r^{n-1}}{\e\mu}(\ti p(v_++\eta^{E}+\bar \eta)-\ti p(v_++\eta^{E}))
%+\frac{1}{\e\mu}(\ti p(v_++\bar \eta++\eta^{E}(1))-\ti p(v_++\eta^{E}(1)))\\
%&-\frac{\e}{\mu}\frac{\chi}{r^{n-1}}+\frac{\e}{\mu}(1-\frac{1}{r^{n-1}})\bar\eta+\frac{\e (n-1)r^{n-1}}{\mu}\int_r^{\infty}\frac{\chi(s)}{s^{2n-1}}ds+\frac{\e (n-1)r^{n-1}}{\mu}\int_r^{\infty}\frac{\bar\eta(s)}{s^{2n-1}}ds=0.
%\end{aligned}
%\end{align}
%\begin{align}
%\label{chinflow}
%\begin{aligned}
%&\chi_r-\frac{r^{n-1}}{\e\mu}(\ti p(v_++\eta^{E}+\bar \eta+\chi)-\ti p(v_++\eta^{E}+\bar \eta)) -\frac{\e}{\mu}\frac{\chi}{r^{n-1}}\\
%&=-\eta^{E}_r+\frac{r^{n-1}}{\e\mu}(\ti p(v_++\eta^{E}+\bar \eta)-\ti p(v_++\eta^{E}))-\frac{1}{\e\mu}(\ti p(v_++\bar \eta++\eta^{E}(1))-\ti p(v_++\eta^{E}(1)))\\
%&-\frac{\e}{\mu}(1-\frac{1}{r^{n-1}})\bar\eta-\frac{\e (n-1)r^{n-1}}{\mu}\int_r^{\infty}\frac{\chi(s)+\bar\eta(s)}{s^{2n-1}}ds
%\end{aligned}
%\end{align}

\begin{align}
\label{inf0bl}
 \left\{
\begin{aligned}
&\chi_r+\frac{r^{n-1}}{\e\mu}a(r)\chi=-\eta^{E}_r+F_1+F_2,\quad r \ge 1,\\[2mm]
&\displaystyle{\lim_{r\to \infty}}\chi(r) =0,\quad \chi(1)=0.
\end{aligned}
 \right.
\end{align}
where 
\begin{align*}
\begin{aligned}
&a := -\frac{(\ti p(v_++\eta^{E}+\bar \eta+\chi)-\ti p(v_++\eta^{E}+\bar \eta))}{\chi}-\frac{\e^2}{r^{2(n-1)}},\\
& F_1:=\frac{r^{n-1}}{\e\mu}(\ti p(v_++\eta^{E}+\bar \eta)-\ti p(v_++\eta^{E}))\\
&\qquad -\frac{1}{\e\mu}(\ti p(v_++\bar \eta+\eta^{E}(1))-\ti p(v_++\eta^{E}(1)))
-\frac{\e}{\mu}(1-\frac{1}{r^{n-1}})\bar\eta,\\
&F_2:=-\frac{\e (n-1)r^{n-1}}{\mu}\int_r^{\infty}\frac{(\eta-\eta^{E})(s)}{s^{2n-1}}ds.
\end{aligned}
\end{align*}
By the same arguments as in the previous sections, 
%$\tilde p'(v)<0\ (v>0)$, Theorem 2.1, and (\ref{4decaybl}), 
it is easy to see that there exist a positive constant $\delta$ which is independent of
$\mu$ satisfying
\begin{equation}
\label{|a|inflow}
a(r) \ge \delta, \qquad r\ge 1,
\end{equation}
for suitably small $|v_--v_+|$ and $|u_-|$.
The theorem for the problem (\ref{inf0bl}) which we need to prove is
\begin{thm} 
\label{ivli}
Let $n \ge 2$ and $\al \in (0,1)$. Then, for any $v_+>0$, there exist
positive constants $\epsilon_0$ and $C$ which are independent of $\mu\in (0,1]$
such that if $|u_-|+|v_--v_+| \le \epsilon_0$, 
the solution $\chi$  of {\rm (\ref{inf0bl})} satisfies
\begin{align}
\label{invsidlimitinf}
\sup_{r\in [1,1+\mu^{\alpha}]}|\chi(r)| \le C\mu^\al.
\end{align}
\end{thm}
\noindent
%We see that Main Theorem \ref{bl+E} is a direct consequence of the Theorem \ref{ivli}. \\

{\it Proof }\quad Applying the Duhamel's principle to the equation (\ref{inf0bl}) 
in terms of $\chi$, we obtain 
\begin{align}
\label{chiduabl}
\chi(r)=\int_1^r e^{-\int_{\tau}^r \frac{a(s)}{\e\mu}s^{n-1}ds}
(-\eta^{E}_{r}+F_1+F_2)(\tau)\,d\tau.
\end{align}
By the same calculation as in (\ref{I2bleu}), 
we can estimate the first term on the right hand side of (\ref{chiduabl}) as 
\begin{align}
\label{4first}
\begin{aligned}
&|\int_1^r e^{-\int_{\tau}^r \frac{a(s)}{\e\mu}s^{n-1}ds}\eta^{E}_r({\tau})\,d\tau|
\le \frac{\e\mu}{\delta}\sup_{r\ge 1}|\eta^{E}_r(r)| \le C\mu.
\\
\end{aligned}
\end{align}
Noting that 
\begin{align*}
\begin{aligned}
|F_1|&\le |\frac{1}{\e\mu}(\ti p(v_++\bar \eta+\eta^{E})-\ti p(v_++\bar \eta+\eta^{E}(1)))|\\
&\quad+ |\frac{1}{\e\mu}(\ti p(v_++\eta^{E})-\ti p(v_++\eta^{E}(1)))|\\
&\quad +|(\frac{r^{n-1}}{\e\mu}-\frac{1}{\e\mu})(\ti p(v_++\bar \eta+\eta^{E})
-\ti p(v_++\eta^{E}))|\\
&\quad +|\frac{\e}{\mu r^{n-1}}(r^{n-1}-1)\bar\eta|\\
%&\quad \le \frac{C}{\mu\e}|r-1|+\frac{r^{n-1}-1}
%{\mu\e}\left(C+\frac{\e^2}{r^{n-1}}\right)\bar\eta 
&\le \frac{C}{\mu\e}|r-1|,
\end{aligned}
\end{align*}
the second term on the right hand side of (\ref{chiduabl}) can be estimated as 
\begin{align}
\label{blF1}
\begin{aligned}
&|\int_1^r e^{-\int_{\tau}^r \frac{a(s)}{\e\mu}s^{n-1}ds}F_1(\tau)\,d\tau| 
\le \frac{C}{\mu\e}\int_1^r e^{-\frac{\delta}{\e\mu}(r-\tau)}|\tau-1|\,d\tau 
\le C\mu^{\alpha},
\end{aligned}
\end{align}
where we used $r\in [1, 1+\mu^{\al}]$. 
By integration by parts, we can estimate the third term of (\ref{chiduabl}) as 
\begin{align}
\label{blF2-1}
\begin{aligned}
&|\int_1^r e^{-\int_{\tau}^r \frac{a(s)}{\e\mu}s^{n-1}ds}F_2(\tau)\,d\tau| \\
&\le \frac{\e (n-1)}{\mu}\int_1^r e^{-\frac{\delta (r^n-\tau^n)}{\e\mu n}}\tau^{n-1}
\int_{\tau}^{\infty}\frac{|(\eta-\eta^{E})(s)|}{s^{2n-1}}ds \, d\tau \\
&\le \frac{\e (n-1)}{\mu}\int_1^r \frac{\e\mu}{\delta}
(e^{-\frac{\delta (r^n-\tau^n)}{\e\mu n}})_{\tau}
\int_{\tau}^{\infty}\frac{|(\eta-\eta^{E})(s)|}{s^{2n-1}}ds \, d\tau \\
&\le \frac{\e^2 (n-1)}{\delta}\left\{\int_{r}^{\infty}\frac{|(\eta-\eta^{E})(s)|}{s^{2n-1}}ds
+\int_1^r  e^{-\frac{\delta (r^n-\tau^n)}{\e\mu n}}
\frac{|(\eta-\eta^{E})(\tau)|}{\tau^{2n-1}} d\tau \right\} \\
&\le \frac{\e^2 (n-1)}{\delta}\left\{\int_{r}^{\infty}\frac{|(\eta-\eta^{E})(s)|}{s^{2n-1}}ds
+\frac{\e\mu}{\delta} \sup_{r\ge 1}|(\eta-\eta^{E})(r)|\right\}\\
&\le C\int_{r}^{\infty}\frac{|(\eta-\eta^{E})(s)|}{s^{2n-1}}ds + C\mu.
\end{aligned}
\end{align}
Further more, the first term on the last line of (\ref{blF2-1}) can be estimated as 
\begin{align}
\label{blF2-2}
\begin{aligned}
&\int_{r}^{\infty}\frac{|(\eta-\eta^{E})(s)|}{s^{2n-1}}\,ds
=\int_{r}^{1+\mu^{\al}}\frac{|(\eta-\eta^{E})(s)|}{s^{2n-1}}\,ds
+\int_{1+\mu^{\al}}^{\infty}\frac{|(\eta-\eta^{E})(s)|}{s^{2n-1}}\,ds\\
%&\le \int_{1}^{1+\mu^{\al}}\frac{|\chi(s)+\bar\eta(s)|}{s^{2n-1}}ds
%+\int_{1}^{\infty}\frac{|\chi(s)+\bar\eta(s)|}{s^{2n-1}}ds\\
&\le \sup_{r\ge 1}|(\eta-\eta^{E})(r)|\int_{1}^{1+\mu^{\al}}1 \,ds
+\sup_{r \ge 1+\mu^{\al}}|(\eta-\eta^{E})(r)|\int_{1}^{\infty}\frac{1}{s^{2n-1}}\,ds\\[3mm]
&\le C\mu^\al +C\mu \le C\mu^\al,
\end{aligned}
\end{align}
where we used the result in Theorem 4.3. 
%where we use the definition $\chi(r):=\eta(r)-(\eta^{E}(r)+\bar \eta(r))$. 
Substituting (\ref{4first}), (\ref{blF1}), (\ref{blF2-1}) and (\ref{blF2-2}) 
into (\ref{chiduabl}), we have
\begin{align}
\label{limitonb}
\begin{aligned}
&\sup_{r\ge 1}|\chi(r)|
\le C\mu^\al.
%&\le \frac{C\mu^{\alpha}}{\delta}+\frac{\e^2 (n-1)}{\delta}\left((\mu^{\al}+\frac{C\e\mu}{\delta})\sup_{1\le r\le 1+\mu^{\al}}|\chi(r)+\bar\eta(r)|
%+C\sup_{1+\mu^{\al}\le r}|\eta(r)-\eta^{E}(r)|\right),
\end{aligned}
\end{align}
Thus, the proof of Theorem 4.4 is completed. \ $\Box$

\bigskip

Now, recalling the fact
$$
|\bar \eta(r)|\le C\exp\{-\nu_0\mu^{\al-1}/\e\},\quad r\ge 1+\mu^\al,
$$
and combining the results in Theorem 4.3 and 4.4, we easily have
%Finally taking limit of $\mu$ to $0$, the right hand side of (\ref{limitonb}) tends to $0$, that is 
%\begin{align}
%\label{limitchil}
%\begin{aligned}
%\sup_{r\ge 1+\mu^{\alpha}}|\chi(r)| \to 0, \quad (\mu \to 0),
%\end{aligned}
%\end{align}
%$\Box$\\
%(\ref{limitchir}) and (\ref{limitchil}) means that 
%which means that 
%\begin{align*}
%\label{limitonr}
%\begin{aligned}
%\sup_{1+\mu^{\alpha}\le r}|\eta(r)-\eta^{E}(r)| \to 0, \quad (\mu \to 0).
%\end{aligned}
%\end{align*}
%and
\begin{align*}
\begin{aligned}
\sup_{r\ge 1}|\eta(r)-\eta^{E}(r)-\bar \eta(r)| \le C\mu^\al,
\end{aligned}
\end{align*}
that is,
\begin{align}
\label{limitonr}
\begin{aligned}
\sup_{r\ge 1}|\eta(r)-\eta^{E}(r)- {\eta}^B(\frac{r-1}{\mu}) 
+ \eta^{E}(1)| \le C\mu^\al.
\end{aligned}
\end{align}
Thus, recalling again
$$
{\rho}(r) = \frac{1}{v_++{\eta}(r)},\quad
{\rho}^{E}(r) = \frac{1}{v_++{\eta}^{E}(r)},\quad
{\rho}^B(y) = \frac{1}{v_++\eta^B(y)},
$$
and
$$
u(r)= \frac{u_-(v_++{\eta}(r))}{v_-r^{n-1}},\quad
{u}^{E}(r) = \frac{u_-(v_++{\eta}^{E}(r))}{r^{n-1}(v_++\eta^{E}(1)},\quad
{u}^B(y) = \frac{u_-(v_++\eta^B(y))}{v_-},
$$
we can complete the proof of the statement (III) in Theorem \ref{bl+E}.
We omit the details.

%%%%%%%%%%%%%%%%%%%%%%%%%%%%%%%%%%%%%%%%%%%%%%%%%%%%%%%%%%%%%%%%%%%

\bigskip

\end{document}